\documentclass[12pt]{article}

\usepackage[latin1]{inputenc}
\usepackage{amsfonts,epsfig,graphicx}
\usepackage{amsmath,amssymb,amsthm}
\usepackage{epsf}
\usepackage{graphics}
\usepackage{amsfonts}
\usepackage{amsmath}
\usepackage{amssymb}
\usepackage[usenames,dvipsnames]{color}
\usepackage{graphicx}
\usepackage{times}

\usepackage{latexsym}

\usepackage{mathbbol}
\usepackage{bbm}
\newtheorem{lemma}{Lemma}
\newtheorem{theorem}{Theorem}
\newtheorem{definition}{Definition}

\newtheorem{example}{Example}

\newtheorem{remark}{Remark}

\usepackage[colorlinks,citecolor=blue]{hyperref}
\usepackage{tikz}
\usetikzlibrary{calc,positioning}
\usepackage{bbm}
\usepackage{float}
\usepackage{enumitem}
\usepackage{comment}

\usepackage{MnSymbol}

\long\def\ignore#1{}
\newcommand{\bdfn}{\begin{definition}}
\newcommand{\edfn}{\end{definition}}

\newcommand{\beq}{\begin{equation}}
\newcommand{\eeq}{\end{equation}}
\newcommand{\brmk}{\begin{remark}}
\newcommand{\ermk}{\end{remark}}
\newcommand{\bthm}{\begin{theorem}}
\newcommand{\ethm}{\end{theorem}}

\newcommand{\epf}{ $~$ \hfill $\blacksquare$}
\newcommand{\bex}{\begin{example}}
\newcommand{\eex}{\end{example}}
\newcommand{\blem}{\begin{lemma}}
\newcommand{\elem}{\end{lemma}}

\newcommand{\mcX}{\mathcal{X}}
\newcommand{\mcA}{\mathcal{A}}
\newcommand{\mcB}{\mathcal{B}}
\newcommand{\mcU}{\mathcal{U}}
\newcommand{\mcM}{\mathcal{M}}
\newcommand{\mcE}{\mathcal{E}}

\newcommand{\mbbR}{\mathbb{R}}
\newcommand{\mbbone}{\mathbb{1}}

\newcommand{\tilX}{\tilde{X}}
\newcommand{\tilu}{\tilde{u}}

\newcommand{\parrzero}{{(0)}}

\newcommand{\parrk}{{(k)}}
\newcommand{\parra}{{(a)}}
\newcommand{\parrb}{{(b)}}

\newcommand{\sqrrk}{{[k]}}

\newcommand{\ringB}{\mathring{B}}

\newcommand{\barmu}{\bar{\mu}}
\newcommand{\bareta}{\bar{\eta}}

\usepackage{tikz,pgfplots}
\usepgfplotslibrary{polar}
\usepgflibrary{shapes.geometric}
\usetikzlibrary{calc}

\pgfplotsset{my style/.append style={axis x line=middle, axis y line=middle, xlabel={$x$}, ylabel={$y$}, axis equal }}

\pgfplotsset{graphstyle01/.append style={axis x line=middle, axis y line=middle, xlabel={$x$}, ylabel={$y$}, axis equal }}

\begin{document}

\begin{center}
{\Large \bf Reversible Markov decision processes and the Gaussian free field}\\
~\\
Venkat Anantharam\\
{\sl EECS Department}\\
{\sl University of California}\\
{\sl Berkeley, CA 94720, U.S.A.}\\
~\\
{\em (Dedicated to the memory of Aristotle (Ari) Arapostathis)}
~\\
~\\
{\bf Abstract}
\end{center}

A Markov decision problem is called reversible if the
stationary controlled Markov chain is reversible under every stationary Markovian strategy. A natural application in which such problems arise
is in the control of Metropolis-Hastings type dynamics. We
characterize all discrete time reversible Markov decision processes with finite state and actions spaces.
We show that policy iteration algorithm for finding an optimal policy
can be significantly simplified Markov decision problems of this type.
We also highlight the relation between the finite time evolution of 
the accrual of reward and the Gaussian free field associated to the controlled Markov chain.

%\maketitle

\iffalse
\begin{flushleft}
{\bf Reversible Markov decision processes and the Gaussian free field}\\
{\sl Venkat Anantharam}\\
{\em 12 May 2022}\\
\end{flushleft}
\fi

\iffalse
\begin{flushleft}
{\em Dedicated to the memory of Aristotle (Ari) Arapostathis}
\end{flushleft}
\fi

\section{Introduction}

We study Markov decision processes (MDPs) in a finite-state, finite-action
%discrete-time 
framework with an average-reward criterion, when the controlled Markov chain is irreducible and reversible in stationarity under every stationary Markov control strategy. This problem was originally studied 
%quite thoroughly 
in special cases
by Cogill and Peng \cite{CP2013}, 
but that work does not seem to have attracted much attention.
%perhaps because of the rather restrictive nature of the assumptions. 
%Our main contribution is to 
%go beyond the work in \cite{CP2013} by clarifying some of the issues left unresolved there
%- in particular, we 
We strengthen the main theorems in \cite{CP2013} by getting rid of superfluous assumptions.
We characterize the class of all such problems.
We also 
%make some interesting 
highlight the
connections between such problems and
%between the control problem and 
%with 
the Gaussian free field 
%associated to a canonical 
of a weighted graph.
%irreducible transient continuous-time Markov chain that can be associated to any such control problem. 

This paper is dedicated to the memory of 
Ari Arapostathis, a good personal friend, who was fond both
of discrete-state MDPs
%in both the discrete world
and of the control problems arising in the Gaussian world of diffusion processes. We hope that the mix of 
%finite-state, finite-action Markov decision theory 
MDPs
with Gaussianity appearing 
in this paper
%in this paper 
%in this class of control problems
-- which is of a form that is unusual in the control context
%although the connection with Gaussianity is
%quite well-known to probabilists studying reversible Markov processes -- 
--
would have met with his approval.

\section{Setup}

$\mcX$ and $\mcU$ are finite sets, denoting the set of states and the set of actions respectively. 
To avoid 
%having to repeatedly deal 
dealing
with corner cases, we assume that both $\mcX$ and $\mcU$ have cardinality at least $2$.
For each $u \in \mcU$ 
%we are given a 
%transition probability matrix 
let
$P(u) := \left[ \begin{array}{c}  p_{ij}(u) \end{array}  \right]$ 
be a transition probability matrix (TPM) on $\mcX$, where $p_{ij}(u)$ denotes the conditional probability that the next state is $j$ when the current state is $i$ and the action taken is $u$. 
We are also given a function $r: \mcX \times \mcU \to \mbbR$, where $r(i,u)$ denotes the reward received if the
current state is $i$ and the current action is $u$.

A {\em stationary randomized Markov strategy} $\mu$ is defined to be a choice of conditional probability distributions
$(\mu(u|i):  u \in \mcU, i \in \mcX)$ and 
results in the 
TPM
%transition probability matrix
$P(\mu) := \left[ \begin{array}{c} p_{ij}(\mu) \end{array}
 \right]$ on $\mcX$, where 
\[
p_{ij}(\mu) := \sum_u p_{ij}(u) \mu(u|i).
\]
We write $\mcM$ for the set of stationary randomized Markov control strategies.
The interpretation of $P(\mu)$ is as the 
TPM
%transition probability matrix 
of the controlled Markov chain when the 
%stationary Markov 
strategy $\mu \in \mcM$ is implemented.
%i.e. when the control strategy is that when in state $i$ the action $u$ is chosen with probability $\mu(u|i)$.

We make the assumption that for each $\mu \in \mcM$ the 
TPM
%transition probability matrix 
$P(\mu)$ is irreducible and reversible. We will then say that we are dealing with a {\em reversible Markov decision problem} (RMDP).
This assumption may seem quite restrictive, but 
%as we will see shortly the class of scenarios covered
%- which we can characterize completely -
it 
seems to be sufficiently interesting to 
merit some attention.
For instance, applications to the optimal design of algorithms of the Metropolis-Hastings type to generate a target probability distribution on a large set of combinatorial configurations, i.e. the Markov Chain Monte Carlo method, are discussed in some depth in \cite[Sec. 5]{CP2013}.

By the assumption of irreducibility there is unique probability distribution 
$ \pi(\mu) := (\pi_i(\mu):  i \in \mcX)$, called the {\em stationary distribution} of $P(\mu)$,
which can be thought of as a column vector satisfying $\pi(\mu)^T P(\mu) = \pi(\mu)^T$.
The assumption of reversibility says that we have 
\begin{equation}        \label{eq:detbal}
\pi_i(\mu) p_{ij}(\mu) = \pi_j(\mu) p_{ji}(\mu), \mbox{ for all $i,j \in \mcX$}.
\end{equation}
The conditions in \eqref{eq:detbal} are often called a {\em detailed-balance} assumption. 
Note that $(\pi_i(\mu) p_{ij}(\mu): (i,j) \in \mcX \times \mcX))$ is a probability distribution,
called the {\em occupation measure} of $P(\mu)$, and the reversibility assumption for $P(\mu)$ is 
equivalent to the assumption that the occupation measure is symmetric when viewed as a matrix.

We will denote the set of stationary deterministic Markov control strategies by $\bar{\mcM}$. 
%We write $\bar{\mcM}$ for the set of stationary deterministic
%Markov control strategies. 
and write $\barmu$ for such a strategy.
Thus $\barmu \in \bar{\mcM}$ is a function $\barmu : \mcX \to \mcU$ and,
with an abuse of notation,
can also be thought of as 
the stationary randomized Markov control strategy $\barmu$ where 
$\barmu(u|i)$ equals $1$ if $u = \barmu(i)$ and $0$ otherwise.
%$\barmu(u|i) = \delta_{\barmu(i),u}$,
%with $\delta_{a,b}$ for $a$ and $b$ in some fixed finite set $\mcA$ denoting the Kronecker delta function, given
%by $\delta_{a,b} = 1$ if $a = b$ and $\delta_{a,b} = 0$ otherwise.
Note that $|\bar{\mcM}| = | \mcU |^{| \mcX |}$, where $|\mcA|$ denotes the cardinality of a finite set
$\mcA$. Of course, $P(\barmu)$ need not be distinct for distinct $\barmu \in \bar{\mcM}$. 
Similarly, $\mcM$ can be thought of as the product of $\mcX$ copies of the probability simplex based on $\mcU$.

%brmk

Even though irreducibility of the $P(u)$ for $u \in \mcU$ is not explicitly mentioned as a condition 
in the definition of the notion of 
%a reversible Markov decision problem 
an RMDP
in \cite{CP2013}, it seems to be 
implicitly assumed, since the notion of reversibility seems to be discussed
there under the implicit assumption that
there is a unique stationary distribution. Thus the use of the terminology 
``reversible Markov decision process"
in this document seems to be consistent with its use in \cite{CP2013}.

%\ermk

%In Section xxx we give some concrete instances of 
%reversible Markov decision problems 
%RMDPs
%of the type described in Example \ref{ex:example} 
%to give some feeling for the potential relevance of this class of problems.
%\color{red}
%Need to write this section.
%\color{black}

\section{Initial results}

Our first claim is the following simple observation. For completeness, a proof is provided in Appendix
\ref{app:detvsrand}.

\blem       \label{lem:detvsrand}

$P(\mu)$ is irreducible and reversible for each $\mu \in \mcM$
iff $P(\barmu)$ is irreducible and reversible for each $\barmu \in \bar{\mcM}$.
\epf

\elem

As pointed out in \cite{CP2013}, a natural class of examples of 
RMDPs
%Markov decision problems 
%satisfying the irreducibility and reversibility conditions above 
arises as follows.

\bex        \label{ex:example}

Let $G := (\mcX, \mcE)$ be a simple connected graph 
%with no self-loops 
with the finite vertex set $\mcX$ and edge set 
$\mcE$. (Recall that a graph is called simple if it does not have multiple edges between any pair of vertices
and does not have any self-loops.)
To each edge $(i,j) \in \mcE$ (between the vertices $i,j \in \mcX$) associate the 
strictly positive weight $s_{ij}$ (thus $s_{ij} = s_{ji}$). 
%We assume that $G$
%has no self-loops, and so $s_{ii} = 0$ for all $i \in \mcX$.
Since $G$
has no self-loops, we have $s_{ii} = 0$ for all $i \in \mcX$.
Write $s_i$ for $\sum_{j \in \mcX} s_{ij}$ and $S$ for $\sum_{i \in \mcX} s_i$.
Let $P^\parrzero$ denote the 
transition probability matrix on $\mcX$ with 
\[
p_{ij}^\parrzero = \frac{s_{ij}}{s_i}, \mbox{ for all $i, j \in \mcX$}.
\]

Let $\rho: \mcX \times \mcU \to (0,1]$ be given. When the control action is $u \in \mcU$, assume that the
state transitions occur according to $P(u)$, where
\begin{eqnarray*}
p_{ii}(u) &=& 1 - \rho(i,u),\\
p_{ij}(u) &=& \rho(i,u) p_{ij}^\parrzero, \mbox{ if $j \neq i$}.
\end{eqnarray*}

Finally, assume that a reward function $r: \mcX \times \mcU \to \mbbR$ is given.

To check that 
%the conditions are satisfied by such a Markov decision problem,
this results in 
%a reversible Markov decision problem
an RMDP
we first observe that
$P^\parrzero$ is an irreducible and reversible 
TPM
%transition probability matrix 
on $\mcX$. The irreducibility 
is obvious. Reversibility can be checked by observing that the stationary distribution of 
$P^\parrzero$, i.e. $(\pi_i^\parrzero: i \in \mcX)$, is given by $\pi_i^\parrzero = \frac{s_i}{S}$ for $i \in \mcX$.

Given $\mu \in \mcM$, write $\rho(i,\mu)$ for $\sum_{u \in \mcU} \rho(i,u) \mu(u|i)$. for $j \neq i$,
Then it can be checked that we have $p_{ij}(\mu) = \rho(i,\mu) p_{ij}^\parrzero$ for $j \neq i$,
while $p_{ii}(\mu) = 1 - \rho(i,\mu)$.
Now, for each $\mu \in \mcM$, $P(\mu)$ is irreducible since, by assumption,
we have $\rho(i,u) > 0$ for all $(i,u) \in \mcX \times \mcU$.
To check that $P(\mu)$ is reversible it suffices to observe that its stationary distribution,
i.e. $(\pi_i(\mu): i \in \mcX)$, is
given by $(K(\mu) \frac{\pi_i^\parrzero}{\rho(i,\mu)}: i \in \mcX)$, where 
$K(\mu) := \left( \sum_{i \in \mcX} \frac{\pi_i^\parrzero}{\rho(i,\mu)} \right)^{-1}$ is the normalizing constant.
\epf

\eex

%\brmk

In the scenario of Example \ref{ex:example}, if one scales all the weights $s_{ij}$ by the same
positive constant then, with the same $\rho: \mcX \times \mcU \to (0,1]$ and
$r: \mcX \times \mcU \to \mbbR$, one gets the same 
%reversible Markov decision problem, 
RMDP,
since all the $s_i$ and $S$ also scale by the same constant. 
%Indeed, what matters is not so much 
%the weighted graph but the 
What matters is the irreducible reversible transition probability matrix $P^\parrzero$ with zero 
diagonal entries 
defined by the weighted graph.
%which it defines.
% Indeed, it is straightforward to 
Conversely, one can
 check that any irreducible
reversible 
TPM
%transition probability matrix 
$P^\parrzero$ with entries $p_{ij}^\parrzero$, $i,j \in \mcX$ and zero diagonal entries
can be though of as arising from the simple connected graph $G := (\mcX, \mcE)$ with $(i,j) \in \mcE$
iff $p_{ij}^\parrzero > 0$, with weight $s_{ij} := \pi_i^\parrzero p_{ij}^\parrzero$, where
$(\pi_i^\parrzero: i \in \mcX)$ is the stationary distribution of $P^\parrzero$.
%\epf

%\ermk

%It turns out that one can 
%reversible Markov decision problem. 
%RMDP.
As stated in the following simple
lemma, whose proof is in Appendix \ref{app:thereisagraph},
one can associate a simple connected graph 
to any RMDP.
We will refer to this graph as the 
{\em canonical graph} of the RMDP.

\blem       \label{lem:thereisagraph}

Consider 
%a reversible Markov decision problem,
an RMDP,
defined by $(P(u): u \in \mcU)$ and $r: \mcX \times \mcU \to \mbbR$ as 
above. Then there must exist a simple connected graph $G := (\mcX, \mcE)$
%without self-loops
such that for all $u \in \mcU$ and distinct $i,j \in \mcX$
we have $p_{ij}(u) > 0$ iff $(i,j) \in \mcE$.
\epf

\elem

In fact, as stated in the following theorem, it turns out that under relatively mild conditions 
every 
%reversible Markov decision problem 
RMDP
must be of the 
form described in Example \ref{ex:example}. The proof is provided in Appendix \ref{app:biconnected}.
%To understand the statement of the theorem we first recall the following standard definition.

\bthm       \label{thm:biconnected}

Consider 
%a reversible Markov decision problem,
an RMDP,
defined by $(P(u): u \in \mcU)$ and $r: \mcX \times \mcU \to \mbbR$.
%Then there must exist a simple connected graph 
Let $G := (\mcX, \mcE)$
%without self-loops
%such that for all $u \in \mcU$ and distinct $i,j \in \mcX$
%we have $p_{ij}(u) > 0$ iff $(i,j) \in \mcE$.
be the canonical 
%simple connected graph 
graph
associated to this problem, as in 
Lemma \ref{lem:thereisagraph}.
Suppose now that this graph is {\em biconnected}.
(Recall that a graph is called biconnected -- or {\em $2$-connected} -- if whenever any 
single vertex, together with all the edges involving that vertex, is removed
the resulting graph continues to be connected.)
Then 
%one can find
%strictly positive weights $s_{ij}$ for $(i,j) \in \mcE$,
there is an irreducible reversible 
TPM
%transition probability matrix
$P^\parrzero$ on $\mcX$ such that $p_{ij}^\parrzero > 0$ iff $(i,j) \in \mcE$,
and $\rho: \mcX \times \mcU \to (0,1]$, such that 
for each $u \in \mcU$ we have
$p_{ij}(u) = \rho(i,u) p_{ij}^\parrzero$ for $j \neq i$,
and $p_{ii}(u) = 1 - \rho(i,u)$.
%where
%$p_{ij}(0) := \frac{s_{ij}}{s_i}$ for $j \neq i$, with
%$s_i := \sum_{j \in \mcX} s_{ij}$ for $i \in \mcX$.
\epf
\ethm

%\brmk

Much of the discussion in \cite{CP2013} centers around 
%reversible Markov decision problems 
RMDPs which have a Hamilton cycle in their canonical graph.
These are biconnected, and hence
of the kind in 
Example \ref{ex:example}.
%(and, as we will soon see, with superfluous conditions). 
 However, 
as seen from Example \ref{ex:counterexample} below, there are 
%reversible Markov decision problems 
RMDPs
that
are not of the type in Example \ref{ex:example}.
%\color{red}
%May want to delete this remark.
%\color{black}

%\ermk

\bex        \label{ex:counterexample}

%The assumption of $2$-connectedness is essential for Theorem \ref{thm:biconnected}
%to hold, as seen from the following example. 
Let $\mcX = \{1,2,3\}$ and $\mcU = \{1,2\}$.
Choose $a \neq b$ such that $0 < a, b, < 1$. Let
\[
P(1) = \left[ \begin{array}{ccc} 0 & a & 1-a\\ 1 & 0 & 0 \\ 1 & 0 & 0 \end{array} \right]
\mbox{ and } 
P(2) = \left[ \begin{array}{ccc} 0 & b & 1-b\\ 1 & 0 & 0 \\ 1 & 0 & 0 \end{array} \right].
\]
It can be checked that for any 
%stationary randomized Markov strategy 
$\mu \in \mcM$ we have
\[
P(\mu) = \left[ \begin{array}{ccc} 0 & \lambda a + (1-\lambda) b & \lambda (1-a) + (1-\lambda) (1-b) 
\\ 1 & 0 & 0 \\ 1 & 0 & 0 \end{array} \right],
\]
for some $\lambda \in [0,1]$ (depending on $\mu$). 
Since $P(\mu)$ is irreducible and reversible, with the stationary distribution being
\[
\pi(\mu) = \left[ \begin{array}{ccc} \frac{1}{2} & \frac{1}{2} \left( \lambda a + (1-\lambda) b \right) & 
\frac{1}{2} \left( \lambda (1-a) + (1-\lambda) (1-b) \right) \end{array} \right],
\]
this, together with some reward function $r: \mcX \times \mcU \to \mbbR$, defines 
%a reversible Markov decision problem. 
an RMDP.
However this 
%reversible Markov decision problem 
RMDP
is not of the form in 
Example \ref{ex:example} as can be seen, for instance, by noticing that $p_{11}(\mu) = 0$ for all
$\mu$ but $p_{12}(\mu)$ can take on distinct values for distinct $\mu$ (since we assumed that $a \neq b$).

Here the 
%simple connected graph $G$ associated to the 
%reversible Markov decision problem 
canonical graph of the
RMDP
has vertex set 
$\mcX$ and edge set $\mcE = \{ (1,2), (1,3) \}$. Note that 
%$G$ is not $2$-connected.
this graph is not biconnected.
\epf

\eex

\section{A characterization of reversible Markov decision problems}

Roughly speaking, a simple connected graph $G := (\mcX, \mcE)$ is as far from
being biconnected as it can be if there is unique path between every pair of vertices
of the graph, i.e. if the graph is a tree. This is of course not precisely true, since
the graph with two vertices connected with a single edge is both biconnected and a tree
and, more generally, in any tree the removal of a leaf vertex together with the edge 
connected to it leaves behind a connected graph. Nevertheless, this rough intuition suggests that
one should pay special attention to trees. As the following simple result shows, in contrast to the
case considered in Theorem \ref{thm:biconnected}, 
%where the graph associated to the 
%reversible Markov decision process 
%RMDP
%is biconnected and this imposes strong restrictions on the structure of the
%Markov decision problem, in case the graph associated to the 
%reversible Markov decision problem 
when the canonical graph of an
RMDP
is
a tree there are hardly any restrictions on the structure of the decision problem. 
The proof is in Appendix \ref{app:treecase}.

\blem       \label{lem:treecase}

Let $G := (\mcX, \mcE)$ be a tree. Let $(P(u): u \in \mcU)$ be any collection of 
TPMs
%transition probability matrices
on $\mcX$ satisfying the condition that $p_{ij}(u) > 0$ iff $(i,j) \in \mcE$. Then, together with
a reward function $r: \mcX \times \mcU \to \mbbR$, this defines 
%a reversible Markov decision problem.
an RMDP.
\epf

\elem

We now proceed to characterize all 
%reversible Markov decision problems. 
RMDPs.
It turns out that the situations discussed in Theorem \ref{thm:biconnected} and Lemma \ref{lem:treecase} are extreme cases and, in a sense, the general case lies between these two extremes. Underlying this is the well-known
{\em block graph} structure of a
simple connected graph $G := (\mcX, \mcE)$.
Recall that a {\em cutvertex} of $G$ 
is a vertex such that if we remove that vertex and the edges connected to it, the resulting graph is disconnected.
A {\em block} is defined to be a maximal connected subgraph of $G$
that has no cutvertices.
Thus a block $B$ is biconnected; in particular, it is either a subgraph comprised of a single edge (in which case it has two vertices) or has the property that 
given any three distinct vertices $i,j,k \in B$ there is path between $j$ and $k$ in $B$ that does not meet $i$. 
Also, any two blocks $B$ and $B^\prime$ that intersect
do so at a uniquely defined vertex,
called an {\em articulation point} of the 
block graph structure. 
An articulation point will be a cutvertex of
$G$ (but not of $B$ or $B^\prime$, since 
$B$ and $B^\prime$ are blocks and so do not have cutvertices). If there is only one block in the block structure, then there are no articulation points. If there is more than one block then every block has at least one articulation point, but in general may have several articulation points. Every articulation point then lies in at least two blocks, but may in general lie in several blocks.
%Here a block is defined to be a maximal connected subgraph without cutvertices of $G := (\mcX, \mcE)$. We now recall the underlying notions and given 
An illustrative example of the block graph structure is given in Figure
\ref{fig:blockstructure} below; see e.g. \cite[Sec. 3.1]{Diestel} for more details
(we focus on connected graphs, even though the
block graph can be defined more generally).

%\color{red}
%Need to add a figure connected to the block graph structure of a connected graph.
%\color{black}

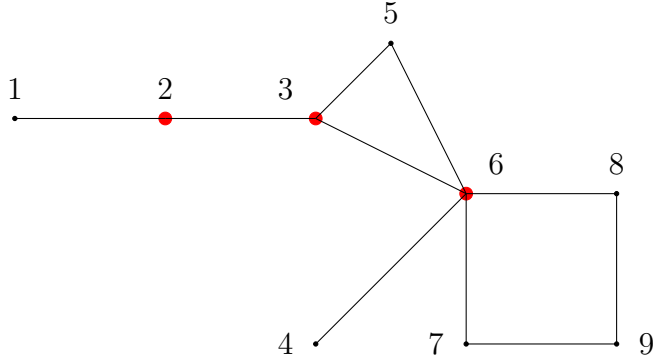
\begin{figure}		\label{fig:blockstructure}
\begin{center}

\begin{tikzpicture}[xscale=4,yscale=4]

%\filldraw [red] (1,1.5) circle (0.5pt);
%\filldraw [red] (-1,1.25) circle (0.5pt);

\iffalse
\filldraw [black] (-1.5,0.5) circle (0.2pt);
\filldraw [black] (-2,0.5) circle (0.2pt);
\filldraw [black] (-1.5,0) circle (0.2pt);
\filldraw [black] (-2,0) circle (0.2pt);
\filldraw [black] (-2.25,1) circle (0.2pt);
\filldraw [black] (-2.5,0.75) circle (0.2pt);
\filldraw [black] (-2.5,0) circle (0.2pt);
\filldraw [black] (-3,0.75) circle (0.2pt);
\filldraw [black] (-3.5,0.75) circle (0.2pt);
\node at (-3.5,0.85) {$1$};
\node at (-3,0.85) {$2$};
\node at (-2.6,0.85) {$3$};
\node at (-2.6,0) {$4$};
\node at (-2.25,1.1) {$5$};
\node at (-1.9,0.6) {$6$};
\node at (-2.1,0) {$7$};
\node at (-1.5,0.6) {$8$};
\node at (-1.4,0) {$9$};
\draw[-] (-1.5,0)--(-2,0);
\draw[-] (-1.5,0.5)--(-2,0.5);
\draw[-] (-1.5,0)--(-1.5,0.5);
\draw[-] (-2,0)--(-2,0.5);
\draw[-] (-2,0.5)--(-2.25,1);
\draw[-] (-2,0.5)--(-2.5,0.75);
\draw[-] (-2.25,1)--(-2.5,0.75);
\draw[-] (-2,0.5)--(-2.5,0);
\draw[-] (-2.5,0.75)--(-3,0.75);
\draw[-] (-3,0.75)--(-3.5,0.75);
\fi

\filldraw [black] (-1.5,-1.5) circle (0.2pt);
\filldraw [red] (-2,-1.5) circle (0.6pt);
\filldraw [black] (-1.5,-2) circle (0.2pt);
\filldraw [black] (-2,-2) circle (0.2pt);
\filldraw [black] (-2.25,-1) circle (0.2pt);
\filldraw [red] (-2.5,-1.25) circle (0.6pt);
\filldraw [black] (-2.5,-2) circle (0.2pt);
\filldraw [red] (-3,-1.25) circle (0.6pt);
\filldraw [black] (-3.5,-1.25) circle (0.2pt);
\node at (-3.5,-1.15) {$1$};
\node at (-3,-1.15) {$2$};
\node at (-2.6,-1.15) {$3$};
\node at (-2.6,-2) {$4$};
\node at (-2.25,-0.9) {$5$};
\node at (-1.9,-1.4) {$6$};
\node at (-2.1,-2) {$7$};
\node at (-1.5,-1.4) {$8$};
\node at (-1.4,-2) {$9$};
\draw[-] (-1.5,-2)--(-2,-2);
\draw[-] (-1.5,-1.5)--(-2,-1.5);
\draw[-] (-1.5,-2)--(-1.5,-1.5);
\draw[-] (-2,-2)--(-2,-1.5);
\draw[-] (-2,-1.5)--(-2.25,-1);
\draw[-] (-2,-1.5)--(-2.5,-1.25);
\draw[-] (-2.25,-1)--(-2.5,-1.25);
\draw[-] (-2,-1.5)--(-2.5,-2);
\draw[-] (-2.5,-1.25)--(-3,-1.25);
\draw[-] (-3,-1.25)--(-3.5,-1.25);

\end{tikzpicture}

\iffalse
\caption{A simple connected graph is in the upper part of the figure. The vertices are numbered as indicated.
The associated block graph is below.
The articulation points are the vertices $2$, $3$, and $6$, and are
depicted by thick red nodes. There are five blocks, namely
$\{1,2\}$, $\{2,3\}$, $\{3,5,6\}$, $\{4,6\}$, and $\{6,7,8,9\}$.
Note that the articulation point $2$ is shared by two blocks, as is the articulation point
$3$, while the articulation point $6$ is shared by three blocks. Note that each 
block has at least one articulation point, while the block $\{3,5,6\}$ has two articulation
points.
 }
 \fi
 \caption{The block graph associated to a simple connected graph with
 nine nodes, numbered as indicated, is shown.
 The articulation points are the vertices $2$, $3$, and $6$, and are
depicted by thick red nodes. There are five blocks, namely
$\{1,2\}$, $\{2,3\}$, $\{3,5,6\}$, $\{4,6\}$, and $\{6,7,8,9\}$.
Note that the articulation point $2$ is shared by two blocks, as is the articulation point
$3$, while the articulation point $6$ is shared by three blocks. Note that
(since there is more than one block) each 
block has at least one articulation point, while the block $\{3,5,6\}$ has two articulation
points.
 }
\end{center}
\end{figure}

Given the simple connected graph 
$G := (\mcX, \mcE)$, we write $\mcA$ for the
set of articulation points and $\mcB$ for the set of blocks. Note that each $a \in \mcA$ is a vertex of $G$, while each $B \in \mcB$ is a subgraph of $G$. Nevertheless, with an abuse of notation, we will 
also use $B$ to denote the vertex set of the block $B$. Thus we write $a \in B$ to indicate that the articulation point $a$ is in the vertex set of $B$ and 
similarly write $B \ni a$ to indicate that the vertex set of $B$ contains the articulation point $a$. 
Further, we will write
$\ringB$ for the subset of those vertices of the block $B$ that are not articulation points.

In the following example we describe a class of 
%reversible Markov decision problems 
RMDPs
that is broader in scope than those considered in Example \ref{ex:example} and Lemma \ref{lem:treecase} (in particular Example \ref{ex:counterexample}), including both of these as special cases.

\bex        \label{ex:genexample}

Let $G := (\mcX, \mcE)$ be 
a simple connected graph. In the block graph structure of $G$, let $\mcA$ denote the set of articulation points and $\mcB$ the set of blocks. For each $B \in \mcB$ let $P^\parrzero(B)$ be a given irreducible reversible
TPM
%transition probability matrix
 on $B$ with $p_{ii}(B) = 0$ for all $i \in B$. For each articulation point $a \in \mcA$ (if any) and $u \in \mcU$, let 
$(\nu_a(u,B): B \ni a)$ be strictly positive numbers satisfying $\sum_{B \ni a} \nu_a(u,B) = 1$. For $i \in \ringB$ and
$u \in \mcU$, define $\nu_i(u,B) = 1$,
and for all $i \in \mcX$ define $\nu_i(u,B) = 0$ 
for all $u \in \mcU$ if $i \notin B$. Let $\rho: \mcX \times \mcU \to (0,1]$ be given.

For each $u \in \mcU$ define 
$P(u) := \left[ \begin{array}{c}  p_{ij}(u) \end{array}  \right]$,
a 
TPM
%transition probability matrix 
on $\mcX$, by 
\begin{eqnarray}    \label{eq:gencaseeqns}
%p_{ij}(u) &=& \rho(i,u) 
%p_{ij}^\parrzero (B)
%\mbox{ if $i \in \ringB$, $j \in B$, $j \neq i$, $B \in \mcB$},\\
%p_{aj}(u) &=& \rho(a,u) \nu_a(u,B) 
%p_{ij}^\parrzero(B), \mbox{ if 
%$a \in B$, $j \in B$, $j \neq a$, $B \in \mcB$},\\
p_{ij}(u) &=& \rho(i,u) \nu_i(u,B)
p_{ij}^\parrzero (B)
\mbox{ if $i \in B$, $j \in B$, $j \neq i$, $B \in \mcB$},\\
p_{ii}(u) &=& 1 - \rho(i,u), 
\mbox{ if $i \in \mcX$}, \nonumber\\
p_{ij}(u) &=& 0, \mbox{ otherwise}. \nonumber
\end{eqnarray}
%Here, for clarity, we have separated out the
%first two lines, though we could have written them together as
%\[
%p_{ij}(u) = \rho(i,u) \nu_i(u,B)
%p_{ij}^\parrzero (B)
%\mbox{ if $i \in B$, $j \in B$, $j \neq i$, $B \in \mcB$}.
%\]
Then, together with
$r: \mcX \times \mcU \to \mbbR$,
this defines 
%a reversible Markov decision problem.
an RMDP.

If there is only a single block, then there are no articulation points and we are in the scenario of Example \ref{ex:example},
where the claim has already been established. Suppose therefore that there are two or more blocks (thus every block has at least one articulation point).
To verify the claim, we need to check that for each $\mu \in \mcM$ the matrix
$P(\mu)$ on $\mcX$ is an irreducible reversible 
TPM.
%transition probability matrix. 
It is straightforward to check that $P(\mu)$ is a 
TPM.
%transition probability matrix.
Noting that for $i \neq j$ we have 
$p_{ij}(\mu) > 0$ iff $(i,j) \in \mcE$, we see that $P(\mu)$ is irreducible.

Let $\tau_i(u,B) := \rho(i,u) \nu_i(u,B)$,
and for $\mu \in \mcM$ let 
$\tau_i(\mu,B) := 
\sum_u \tau_i(u,B) \mu(u|i)$.
Let 
$(\psi_i(B): i \in B)$ denote the stationary probability distribution of $P^\parrzero(B)$. Then, by the assumed reversibility of this matrix, we have
\[
\psi_i(B) p_{ij}^\parrzero(B) = \psi_j(B) p_{ji}^\parrzero(B), \mbox{ for all $i,j \in B$, $B \in \mcB$}. 
\]
Consider the vector 
$(\frac{\psi_i(B)}{\tau_i(\mu, B)}: i \in B)$. Observe now that we have 
\begin{equation}        \label{eq:innerbalance}
\frac{\psi_i(B)}{\tau_i(\mu, B)} p_{ij}(\mu) =
\frac{\psi_j(B)}{\tau_j(\mu, B)} p_{ji}(\mu)
\mbox{ for all $i,j \in B$, $B \in \mcB$}.
\end{equation}
We now claim that we can find positive constants $(m(\mu, B): B \in \mcB)$
such for every $B, B^\prime \in \mcB$,
$B \neq B^\prime$, if they share an articulation point $a \in \mcA$ (i.e. $a \in B$, $a \in B^\prime$), then we have
\begin{equation}        \label{eq:matching}
m(\mu, B) \frac{\psi_a(B)}{\tau_a(\mu, B)}
=
m(\mu, B^\prime) \frac{\psi_a(B^\prime)}{\tau_a(\mu, B^\prime)}.
\end{equation}
Since this number 
does not depend on the choice of $B \in \mcB$ containing $a$, 
%depends only on 
%$a$ (for the given $\mu$), 
let us denote it by $\gamma_a(\mu)$. Let us also write
$\gamma_i(\mu)$ for 
$m(\mu, B) \frac{\psi_i(B)}{\tau_i(\mu, B)}$ for $i \in \ringB$ for any $B \in \mcB$. 
With this notation in place, we further claim that we can choose the 
$(m(\mu, B): B \in \mcB)$ such that 
\begin{equation}        \label{eq:isaprobdist}
\sum_{i \in \mcX} \gamma_i(\mu) = 1. 
\end{equation}
It can then be checked that 
$(\gamma_i(\mu): i \in \mcX)$ is then the stationary distribution of $P(\mu)$ and, based on \eqref{eq:innerbalance}, we can conclude that $P(\mu)$ is reversible.

To find the scaling factors $(m(\mu, B): B \in \mcB)$
with the claimed properties, pick any block and
call it the root. Because there are at least two blocks, this block has at least one articulation point, and each such articulation point is associated with unique block other than the root. Call these blocks the ones at depth $1$. 
If any such block has additional articulation points (other than the one it shares with the root), each of these will be associated with a new block, and we will call the blocks identified in this way (from all the blocks at depth $1$) the blocks at depth $2$, and so on. 
We start with a scaling factor $1$ for the root, and see that we can set the scaling factor uniquely for each of the blocks at depth $1$ in order to get the matching condition in \eqref{eq:matching} to hold at all the articulation points that are shared between the root and the blocks at depth $1$. We can next set the scaling factor for each of the blocks at depth $2$ uniquely in order to get the matching condition in \eqref{eq:matching} to hold at all the articulation points that are shared between a block at depth $1$ and a block at depth $2$ and so on. At the end of this process we have scaling factors such that the condition in \eqref{eq:matching} holds at all articulation points and then we can 
finally scale all the scaling factor jointly
by the same positive constant to get the condition in \eqref{eq:isaprobdist}.
\epf

\eex

%\brmk

The class of 
%reversible Markov decision problems 
RMDPs
arising as in Example \ref{ex:genexample} also includes those arising as in Lemma \ref{lem:treecase}.
This corresponds to the case where every block is a single edge, which is equivalent to the case where the given graph
$G := (\mcX, \mcE)$ is a tree. 
Each $P^\parrzero(B)$ is then of the form
$\left[ \begin{array}{cc} 0 & 1 \\
1 & 0 \end{array} \right]$.
If $|\mcX| =2$ then there are no articulation points and the scenario is covered in Theorem \ref{thm:biconnected} (and also in Lemma \ref{lem:treecase}). If $|\mcX| \ge 3$ the
articulation points are precisely the  non-leaf vertices of the tree.
%\epf

%\ermk

It turns out that the scenarios covered in Example \ref{ex:genexample} completely characterize all the ways in which 
%a reversible Markov decision problem 
an RMDP
can arise. This is stated in the following theorem, whose proof is in Appendix
\ref{app:maintheorem}.

\bthm       \label{thm:maintheorem}

Consider 
%a reversible Markov decision problem,
an RMDP,
defined by the 
TPMs
%transition probability matrices 
$(P(u): u \in \mcU)$ and $r: \mcX \times \mcU \to \mbbR$.
Let $G := (\mcX, \mcE)$
be the canonical 
%simple connected graph 
graph
associated to this problem, as in 
Lemma \ref{lem:thereisagraph}.
In the block graph structure of $G$, let $\mcA$ denote the set of articulation points and $\mcB$ the set of blocks.
Then for each $B \in \mcB$ there will
be an irreducible reversible 
TPM
%transition probability matrix 
$P^\parrzero(B)$
on $B$ with $p_{ii}(B) = 0$ for all $i \in B$; for each articulation point $a \in \mcA$ (if any) and $u \in \mcU$ there will be strictly positive numbers
$(\nu_a(u,B): B \ni a)$,
satisfying $\sum_{B \ni a} \nu_a(u,B) = 1$;
and there will be
$\rho: \mcX \times \mcU \to (0,1]$ 
such that for each $u \in \mcU$ the 
entries of the matrix $P(u)$ are given by 
\eqref{eq:gencaseeqns}.
\epf

\ethm

\section{Dynamic programming equations and policy iteration}       \label{sec:policyiteration}

Consider 
an MDP
%a Markov decision problem (MDP) 
defined by a family $(P(u): u \in \mcU)$ 
where each $P(u)$ is an
irreducible 
TPM
%transition probability matrix 
on $\mcX$, and a reward function
$r: \mcX \times \mcU \to \mbbR$.
Here $\mcX$ and $\mcU$ are finite sets each assumed to be of cardinality at most $2$. 
Given $\mu \in \mcM$, let $\beta(\mu)$ denote the long term average reward associated to the stationary randomized Markov strategy $\mu$.
Then we have $\beta(\mu) = \sum_i \pi_i(\mu) r(i,\mu)$, where $\pi(\mu) = (\pi_i(\mu): i \in \mcX)$ denotes the stationary distribution of $P(\mu)$ and $r(i,\mu) := \sum_u r(i,u) \mu(u|i)$. Further, there is function $h(\mu): \mcX \to \mbbR$ such that for all $i \in \mcX$ we have
\begin{equation}        \label{eq:Poisson}
\beta(\mu) = r(i, \mu(i)) + \sum_j p_{ij}(\mu) \left(h_j(\mu) - h_i(\mu)\right).
\end{equation}
The family of equations \eqref{eq:Poisson},
one for each $i \in \mcX$,
is often viewed as needing to be 
solved 
for $\beta(\mu)$ and $h(\mu) := (h_i(\mu): i \in \mcX)$,
in which case it is called {\em Poisson's equation} associated to the 
TPM
%transition probability matrix 
$P(\mu)$. Note that the number of variables is one more that then number of equations and, indeed, one can add the same fixed constant to each $h_i(\mu)$ in any solution to find another solution. 
%A natural choice for $h(\mu)$ is given by
%\begin{equation}        \label{eq:excessrewards}
%h(\mu) = \lim_{M \to \infty} \frac{1}{M} \sum_{K=1}^M \left( \sum_{k=0}^{K-1} \left( P(\mu)^k - \mbbone \pi(\mu)^T \right) r(\mu) \right),
%\end{equation}
%where $h(\mu)$ and $r(\mu)$ are thought of as column vectors, with
%$r(\mu) = (r(i,\mu(i)): i \in \mcX)$,
%and $\mbbone$ denotes the all-ones column vector. 

A natural choice for $h(\mu)$, thought of as a column vector,
is given by the Ces\`{a}ro limit of the sequence
$( \sum_{k=0}^{K-1} \left( P(\mu)^k - \mbbone \pi(\mu)^T\right)r(\mu), K \ge 1)$,
where $r(\mu)$ is thought of as the column vector with
$r(\mu) = (r(i,\mu(i)): i \in \mcX)$,
and $\mbbone$ denotes the all-ones column vector. 
%The limit on the right hand side of equation \eqref{eq:excessrewards} exists because 
This Ces\`{a}ro limit exists because the sequence
$(\frac{1}{K} \sum_{k=0}^{K-1} P(\mu)^k, K \ge 1)$ converges geometrically fast to $\mbbone \pi(\mu)^T$ as $K \to \infty$.
%The Cesaro average above (and in equation \eqref{eq:excessrewards}) is needed to deal with the phenomenon of periodicity.
Taking the Ces\`{a}ro limit is 
needed 
to deal with
%since the limit may not exist because of 
the phenomenon of periodicity.

The {\em average cost dynamic programming equation} characterizes an optimal stationary
randomized Markov strategy $\mu$ as one having the property that for each $i \in \mcX$ 
if $\mu(u|i) > 0$ then we must have
\begin{equation}    \label{eq:optchar}
r(i, u) + \sum_j p_{ij}(u) \left(h_j(\mu) - h_i(\mu)\right) = 
\max_v \left( r(i, v) + \sum_j p_{ij}(v) \left(h_j(\mu) - h_i(\mu)\right) \right),
\end{equation}
which implies the form in which it is usually written, namely
\begin{equation}        \label{eq:DPeqn}
\beta(\mu) = \max_v \left( r(i, v) + \sum_j p_{ij}(v) \left(h_j(\mu) - h_i(\mu)\right) \right).
\end{equation}
The characterization of optimal stationary randomized Markov strategies in equation \eqref{eq:optchar}
leads to the {\em policy iteration algorithm} to find an optimal stationary deterministic strategy.
Namely, starting with $\barmu^\parrzero \in \bar{\mcM}$, consider the sequence $(\barmu^\parrk \in \bar{\mcM}, k \ge 0)$ where to get $\barmu^{(k+1)}$ from $\barmu^\parrk$ 
we pick some state $i$ (if possible) for which $\text{argmax}_v \left( r(i, v) + \sum_j p_{ij}(v) \left(h_j(\barmu^\parrk) - h_i(\barmu^\parrk)\right) \right)$ does not equal $\barmu^\parrk(i)$, and replace $\barmu^\parrk(i)$ by an action achieving the argmax. It is well-known that we will then have
$\beta_{\barmu^{(k+1)}} > \beta_{\barmu^\parrk}$ (a proof is given in \cite{CP2013}, for instance) 
and that this iteration will terminate in a finite number of steps to a stationary deterministic optimal strategy.

We turn now to the case where the MDP is an RMDP, i.e. when $P(\mu)$ is reversible for all $\mu \in \mcM$.
Consider first the case where the 
canonical 
%simple connected graph
graph
$G := (\mcX, \mcE)$ 
associated to the RMDP is biconnected.
Then, according to Theorem \ref{thm:biconnected}, we have
$\rho: \mcX \times \mcU \to (0,1]$ 
and an irreducible reversible 
TPM
%transition probability matrix 
$P^\parrzero$ on $\mcX$ such that $p_{ij}(u) = \rho(i,u) p_{ij}^\parrzero$ for all distinct $i,j \in \mcX$ and $p_{ii}(u) = 1 - \rho(i,u)$ for all $i \in \mcX$. As stated in the following theorem, the policy iteration algorithm can be dramatically simplified in this case. The proof is in Appendix
\ref{app:policyiterationbic}.

\bthm       \label{thm:policyiterationbic}

Consider an RMDP whose associated canonical %simple connected graph 
graph
is biconnected. Let $P^\parrzero$ and
$\rho: \mcX \times \mcU \to (0,1]$ be 
as in Theorem \ref{thm:biconnected}. Then any sequence 
$(\barmu^\parrk \in \bar{\mcM}, k \ge 0)$
of stationary deterministic Markov strategies, starting from some $\barmu^\parrzero \in \bar{\mcM}$,
where $\barmu^{(k+1)}$ is got from $\barmu^\parrk$ 
by picking some state $i$ (if possible) for which
\begin{equation}        \label{eq:badbic}
\frac{r(i,\barmu^\parrk) - \beta(\barmu^\parrk)}{\rho(i, \barmu^\parrk)} <
\text{argmax}_v \left( 
\frac{r(i,v) - \beta(\barmu^\parrk)}{\rho(i,v)}
\right)
\end{equation}
and replacing $\barmu^\parrk(i)$ by some action achieving the argmax, has the property that 
$\beta(\barmu^{(k+1)}) > \beta(\barmu^\parrk)$,
and this iteration will terminate in a finite number of steps to a stationary deterministic optimal strategy.
\epf

\ethm

%\brmk

A weaker version of Theorem \ref{thm:policyiterationbic} is proved in \cite[Thm. 4.2]{CP2013} under the assumption that there is a Hamilton cycle in the canonical 
%simple connected graph 
graph
associated to the RMDP (which implies, but is a strictly stronger requirement than biconnectedness) and that at each step of the policy iteration the actions at all states are updated simultaneously in a specific way related to this Hamilton cycle, see \cite[Sec. 4.1]{CP2013}.
%\epf

%\ermk

For a general RMDP it turns out that a simplification of policy iteration
similar to that in Theorem \ref{thm:policyiterationbic} is possible
at vertices that are not articulation points. This is stated in the 
following theorem, whose proof is in Appendix 
\ref{app:policyiterationgen}

\bthm		\label{thm:policyiterationgen}

For a general RMDP, 
let $\mcA$ denote the set of articulation points and $\mcB$ the set of blocks in the
 block graph structure of the canonical graph $G := (\mcX, \mcE)$
associated to it.
Let
$(P^\parrzero(B): B \in \mcB)$,
$(\nu_a(u,B): a \in \mcA, B \ni a, u \in \mcU)$,
and
$\rho: \mcX \times \mcU \to (0,1]$ 
be as in Theorem \ref{thm:maintheorem}.
Let $\barmu \in \bar{\mcM}$ and suppose that 
for some $B \in \mcB$ and $k \in \ringB$ we have 
\begin{equation}        \label{eq:badgen}
\frac{r(k,\barmu) - \beta(\barmu)}{\rho(k, \barmu)} <
\text{argmax}_v \left( 
\frac{r(k,v) - \beta(\barmu)}{\rho(k,v)}
\right)
\end{equation}
Let $\barmu^\prime \in \bar{\mcM}$ be defined by setting
$\barmu^\prime(k) = v$ and $\barmu^\prime(j) = \barmu(j)$ for all
$j \neq k$. Then we have $\beta(\barmu^\prime) > \beta(\barmu)$.
\epf

\ethm

%\section{Examples}      %\label{sec:examples}

\section{The Gaussian free field and the generalized second Ray-Knight theorem}     \label{sec:gff}

For every $\mu \in \mcM$ the Ces\`{a}ro limit of the sequence
$( \sum_{k=0}^{K-1} \left( P(\mu)^k - \mbbone \pi(\mu)^T\right), K \ge 1)$
exists and is called the {\em fundamental matrix} associated to $P(\mu)$
\cite[Sec. 2.2.2]{AldousFill}. Denote this matrix by $Z(\mu)$, with
entries $z_{ij}(\mu)$. 
From the discussion in
Section \ref{sec:policyiteration}, note that 
$h(\mu) := Z(\mu) r(\mu)$, together with $\beta(\mu) = \pi(\mu)^T r(\mu)$, solves Poisson's equation for $P(\mu)$, i.e. equation 
\eqref{eq:Poisson}. Thus, understanding the fundamental matrix 
$Z(\mu)$ is central to understanding the dynamics of the RMDP 
under $\mu \in \mcM$.

$Z(\mu)$ is best understood by 
moving to continuous time, replacing the 
TPM
%transition probability matrix
$P(\mu)$ by the rate matrix $P(\mu) - I$, where $I$ denotes
the identity matrix on $\mcX$.  Let $(X_t(\mu), t \ge 0)$ denote
the corresponding continous time Markov chain. Then one can check that
\begin{equation}		\label{eq:firstZformula}
z_{ij}(\mu) = \lim_{T \to \infty} \left( E_i[ \int_0^T 1( X_t(\mu) = j) dt] - \pi_j(\mu) T\right).
\end{equation}

Since $P(\mu)$ is reversible, it is straightforward
to show that the
matrix on $\mcX$ with entries $\frac{z_{ij}(\mu)}{\pi_j(\mu)}$ is symmetric
\cite[Sec. 3.1]{AldousFill}. Based on \eqref{eq:firstZformula},
we may now write, for the choice of $h(\mu)$ above, for each $i \in \mcX$,
the formula
\begin{equation}		\label{eq:rescaled}
h_i(\mu) = \sum_j \frac{z_{ij}(\mu)}{\pi_j(\mu)} \pi_j(\mu) r_j(\mu)
= 
\sum_j  \lim_{T \to \infty} \left(  E_i[ \frac{1}{\pi_j(\mu)} \int_0^T 1( X_t(\mu) = j)dt] - T\right) \pi_j(\mu) r_j(\mu)
\end{equation}

While this may seem a peculiar thing to do, one natural aspect of the
formula on the RHS of \eqref{eq:rescaled} is that
$\pi_j(\mu) r_j(\mu)$ has the interpretation, in continuous time,
of the rate at which reward is generated in stationarity while in state $j$.
Another natural aspect is that the centering of the integral is 
the actual time and not a state-dependent scaled version of it.
However, the real value of this way of writing the formula comes from
the observation that the matrix with entries 
$\frac{z_{ij}(\mu)}{\pi_j(\mu)}$ is a positive semidefinite
matrix  \cite[Eqn. (3.42)]{AldousFill}. 
This means that we can find a 
multivariate mean zero Gaussian random variable, call it 
$(V_i(\mu) : i \in \mcX)$, with this covariance matrix.
This points to an intriguing and unusual connection between 
Gaussianity and Markov decision theory in the case of RMDP.
As we will see shortly, while the $h_i(\mu)$ are expressed as
asymptotic limits in
\eqref{eq:rescaled}, the introduction of Gaussian methods 
gives, in a sense, much more detailed information about the 
behavior of the functions 
$T \to \frac{1}{\pi_j(\mu)} \int_0^T 1( X_t(\mu) = j)$ and
thus a much more detailed picture of the role of the initial condition
in causing deviations from the stationary rate of generation of reward
in an RMDP.

Notice that we have 
$\sum_i \sum_j \pi_i(\mu) \frac{z_{ij}(\mu)}{\pi_j(\mu)} \pi_j(\mu) = 0$,
and so $\sum_i \pi_i(\mu) V_i(\mu) = 0$ as a random variable. 
Thus, to work with $(V_i(\mu) : i \in \mcX)$ involves, in a sense, a choice
of coordinates to capture the underlying multivariate Gaussian structure.
Other natural choices of coordinates are possible. For instance, for each
$k \in \mcX$ we may define the multivariate Gaussian $(V^\sqrrk_i: i \in \mcX)$
via $V^\sqrrk_i : V_i - V_k$ (so the choice of coordinates in this case makes
$V^\sqrrk_k = 0$). 

Instead of making a choice of coordinates, the Gaussian object of interest can 
be constructed in an intrinsic way. One starts with independent mean zero 
Gaussian random variables on the edges of the canonical graph of the RMDP,
with the variance of the Gaussian on edge $(i,j)$ being
$(\pi_i(\mu) p_{ij}(\mu))^{-1}$. To each edge one associates a direction 
in an arbitrary way, with the understanding that traversing the edge along
its direction corresponds to adding this Gaussian, while traversing it in the
opposite direction corresponds to subtracting this Gaussian. One then 
conditions on being in the subspace of $\mbbR^\mcE$ where the total sum of the Gaussians over
every loop in the canonical graph equals zero. This will allow us to construct
a multivariate Gaussian on the vertices of the canonical graph 
with the property that the Gaussian on each edge is the difference
between those at its endpoints. This multivariate Gaussian on the vertices
is defined only up to one degree of freedom and this is what corresponds to the freedom in the choice of coordinates discussed above.
See \cite[Sec. 9.4]{Janson} for more details. This Gaussian object is
called the {\em Gaussian free field} associated to $P(\mu)$. It is 
discussed in many sources, e.g. \cite[Chap. 5]{LeJanSF}, 
\cite[Sec. 2.8]{LPbook}, \cite[Sec. 9.4]{Janson}, 
\cite[Sec. 2.8]{SznitmanETH}.

For each $k \in \mcX$ the representation of the 
Gaussian free field via the multivariate Gaussian $(V^\sqrrk_i: i \in \mcX)$
also has a natural probabilistic interpretation.
Consider the transient continuous time Markov chain on $\mcX$,
with absorption in state $k$, with the rate $\pi_i(\mu) p_{ij}(\mu)$
of jumping from state $i$ to state $j$ for all $i \neq k$. 
Let $g^\sqrrk_{ij}$, for $i, j \neq k$, denote the mean time spent in state $j$ before absorption. Then it can be checked that
the matrix on $\mcX \backslash \{k\}$ with 
entries $g^\sqrrk_{ij}$ is a symmetric positive definite matrix.
It is, indeed, the covariance matrix of $(V^\sqrrk_i: i \neq k)$.
See \cite{LeJanSF} and \cite{SznitmanETH} for more details. 

Let us also observe that the recurrent continuous time Markov
chain $(\tilX_t, t \ge 0)$ on $\mcX$ with the rate of jumping from
state $i$ to state $j$ being $\pi_i(\mu) p_{ij}(\mu)$ for $j \neq i$
satisfies
\begin{equation}		\label{eq:secondZformula}
\frac{z_{ij}(\mu)}{\pi_j(\mu)} = \lim_{T \to \infty} \left( E_i[ \int_0^T 1( \tilX_t(\mu) = j) dt] - T\right).
\end{equation}
This is basically a consequence of \eqref{eq:firstZformula} but is somewhat
more subtle that it might seem. $(X_t, t \ge 0)$ can be coupled to
$(\tilX_t, t \ge 0)$ by creating the latter from the former by stretching
out each duration of time spent in state $i$ by the factor $\pi_i(\mu)^{-1}$,
for each $i \in \mcX$. But then the integral to a fixed time $T$ in 
\eqref{eq:secondZformula} makes the corresponding integral in 
\eqref{eq:firstZformula} be to a random time. Nevertheless, since
we take the asymptotic limit in $T$, \eqref{eq:secondZformula}
follows from \eqref{eq:firstZformula}.

Now, the generalized second Ray-Knight theorem
\cite[Thm. 2.17]{SznitmanETH} gives us the promised insight
into the transient rates at which rewards are generated in the 
individual states. For $i \in \mcX$ and $t \ge 0$, let
$L_{i,t} := \int_0^t 1(\tilX_s = i) ds$. For $k \in \mcX$ and
$s \ge 0$ define
\[
\Gamma_{k,s} := \inf\{ t \ge 0: L_{k,t} \ge s\},
\]
which is the first time at which the time spent in state $k$ by the 
process $(\tilX_t, t \ge 0)$ is at least $s$. 
We then have
\[
\left(L_{i, \Gamma_{k,s}} + \frac{1}{2} (V^\sqrrk_i)^2: i \in \mcX \right) 
\stackrel{d}{=} 
\left( \frac{1}{2} (V^\sqrrk_i + \sqrt{2 s})^2 : i \in \mcX \right)
\]
for all $s  \ge 0$, where $\stackrel{d}{=}$ denotes equality in 
distribution of the vector random variables on each side.
Here $(V^\sqrrk_i: i \in \mcX)$ is the Gaussian free field, as described
earlier, and is assumed to be independent of 
$(L_{i, \Gamma_{k,s}}: i \in \mcX)$, whose law is taken assuming that the process $(\tilX_t, t \ge 0)$ starts at $k \in \mcX$.

This unusual way in which Gaussians plays a role  in the 
context of RMDP to give insight into the transient behavior of
the generation of reward is quite striking. Our purpose in this paper
has only been to highlight this connection. We leave the exploration of its implications 
to future research.

%meaning in terms of the underlying Markov chain. This is best understood by 
%moving to continuous time, replacing the transition probability matrix
%$P(\mu)$ by the rate matrix $P(\mu) - I$, where $I$ denotes
%the identity matrix on $\mcX$.  Let $(X_t(\mu): t \ge 0)$ denote
%the corresponding continous time Markov chain. Then one can check that
%\[
%z_{ij}(\mu) = \lim_{T \to \infty} \left( E_i[ \int_0^T 1( X_t(\mu) = j)] - 
%\pi_j(\mu) T\right).
%\]
%Let $T_k(\mu) := \inf \{ t > 0 : X_t(\mu) = k\}$ denote the first 
%entrance time of state $k$. 
%of 

%\section{Concluding remarks}		\label{sec:concrmk}

%We have developed the theory of reversible Markov decision processes

\section*{Acknowledgements}

This research was supported by NSF grants
CCF-1618145, CCF-1901004, CIF-2007965, and the NSF Science
\& Technology Center grant CCF-0939370 (Science of Information). The author would like to thank Devon Ding for several discussions 
centered around the monographs 
%of Le Jan 
\cite{LeJanSF} and 
%Sznitman 
\cite{SznitmanETH}, and also for reading the completed paper for a sanity check.

%\section*{Appendix}

\appendix

\section{Proof of Lemma \ref{lem:detvsrand}}        \label{app:detvsrand}

%{\em Proof of Lemma \ref{lem:detvsrand}}:

It is well known that the set of occupation measures as $\mu$ ranges over
$\mcM$ is a closed convex set and every extreme
point of this convex set corresponds to the occupation measure of some
$\barmu \in \bar{\mcM}$
(see e.g. \cite[Lemma 5.2]{ABFGM93}).
Now suppose that $P(\barmu)$ is irreducible and reversible for each $\barmu \in \bar{\mcM}$.
Every $\mu \in \mcM$ can be expressed as a finite convex combination 
%$\mu = \sum_{a \in \mcA} \beta^\parra \barmu^\parra$ where $\mcA$ is a finite set,
%$\beta^\parra > 0$ for all $a \in \mcA$, $\sum_{a \in \mcA} \beta^\parra = 1$, and 
%$\barmu^\parra \in \bar{\mcM}$ for all $a \in \mcA$. We then have
%$P(\mu) = \sum_{a \in \mcA} \beta^\parra P(\barmu^\parra)$. 
%of the $\barmu$ as $\barmu$ ranges over $\bar{\mcM}$ and thus
of $(\barmu: \barmu \in \bar{\mcM})$, and thus
$P(\mu)$ is expressed as the corresponding convex combination of $(P(\barmu): \barmu \in \bar{\mcM})$.
It follows that $P(\mu)$
is irreducible. 

But the occupation measure of $P(\mu)$ is in the convex hull of the 
occupation measures $(P(\barmu): \barmu \in \bar{\mcM})$, and these are symmetric by assumption,
so the occupation measure of $P(\mu)$ is symmetric and hence $P(\mu)$ is reversible.
(Note that the convex combination expressing the occupation measure of $P(\mu)$
in terms of the occupation measures of stationary deterministic Markov control strategies may be different
from that used above to express $P(\mu)$ in terms of $(P(\barmu): \barmu \in \bar{\mcM})$.) This proves that $P(\mu)$ 
 irreducible and reversible for each $\mu \in \mcM$. The claim in the opposite direction is 
 obvious. This completes the proof of the lemma.
\epf

\section{Proof of Lemma \ref{lem:thereisagraph}}        \label{app:thereisagraph}

We claim that for every $i \in \mcX$ the 
set of neighbors of $i$ under $(p_{ij}(u): j \in \mcX)$,
namely $\{j \neq i: p_{ij}(u) > 0\}$ is the same for all $u \in \mcU$. 
Suppose, to the contrary, that for some distinct $u, v \in \mcU$ and some
$i \in \mcX$ we have $p_{ij}(u) > 0$ but $p_{ij}(v) = 0$, for some $j \neq i$.
Pick some $\tilu \in \mcU$ (which could be either $u$ or $v$ if desired)
and consider the two stationary deterministic Markov strategies 
$\barmu^\parra$ and $\barmu^\parrb$ given by 
\begin{equation}        \label{eq:twostrats}
\barmu^\parra(i) = u, \barmu^\parrb(i) = v, \mbox{ and } 
\barmu^\parra(l) = \barmu^\parrb(l) = \tilu \mbox{ for $l \neq i$}.
\end{equation}
Since $P(\barmu^\parra)$ is reversible and $p_{ij}(u) > 0$ it follows that
$p_{ji}(\tilu) > 0$. But then, since $P(\barmu^\parrb)$ is reversible, 
it would follow that $p_{ij}(v) > 0$, a contradiction. This establishes the claim.
This also establishes
the existence of a 
simple connected graph $G := (\mcX, \mcE)$
%without self-loops
such that for all $u \in \mcU$ and distinct $i,j \in \mcX$
we have $p_{ij}(u) > 0$ iff $(i,j) \in \mcE$, as claimed.
\epf

\section{Proof of Theorem \ref{thm:biconnected}}        \label{app:biconnected}

Suppose first that all the $P(u)$ for $u \in \mcU$ are the same, and
let $P = \left[ \begin{array}{c} p_{ij} \end{array} \right]$ denote this common
TPM
%transition probability matrix 
over $\mcX$. 
Thus $P$ is irreducible and reversible.
Let $(\pi_i: i \in \mcX)$ denote the stationary distribution of
$P$. We have $\pi_i p_{ij} = \pi_j p_{ji}$ for all $i,j \in \mcX$.

For each $i \in \mcX$ we can choose 
%$\rho(i,u)$ to be the same $\rho(i) \in (0,1]$, where 
%$\rho(i) := 1 - p_{ii}$. 
$\rho(i,u) \in (0,1]$ to be $1 - p_{ii}$.
(Note that we have $p_{ii} < 1$ since $P$ is 
irreducible and $|\mcX| \ge 2$.) 
We then let $p_{ij}^\parrzero := \frac{p_{ij}}{1 - p_{ii}}$ for $i \neq j$, with
$p_{ii}^\parrzero := 0$ for all $i \in \mcX$. It can be checked that this defines an irreducible
TPM
%transition probability matrix 
$P^\parrzero$ on $\mcX$ with $p_{ij}^\parrzero > 0$ iff $(i,j) \in \mcE$,
where $G := (\mcX, \mcE)$ denotes the canonical graph associated to this 
%reversible Markov decision problem. 
RMDP.
It can be checked that the
stationary distribution 
$(\pi_i^\parrzero: i \in \mcX)$ of $P^\parrzero$ is given by $\pi_i^\parrzero = K \pi_i (1 - p_{ii})$, where
$K$ is the normalizing constant. Further, we have 
$\pi_i^\parrzero p_{ij}^\parrzero = \pi_j^\parrzero p_{ji}^\parrzero$ for all $i,j \in \mcX$,
which establishes that $P^\parrzero$ is reversible. 
This completes the proof in this case.
%We then choose $s_{ij} := S \pi_i p_{ij}$
%for $j \neq i$ and $s_{ii} = 0$ for $i \in \mcX$, where 
%$S = \left( \sum_k \pi_k (1- p_{kk}) \right)^{-1}$. 

%From the reversibility of $P$, it can be checked that $s_{ij} = s_{ji}$ for 
%$j \neq i$. Further, we have
%\[
%s_i := \sum_{j \neq i} s_{ij} = S \pi_i (1 - p_{ii}) = \rho(i) S \pi_i,
%\]
%from which it follows that $\rho(i) \frac{s_{ij}}{s_i} = p_{ij}$ for all
%$j \neq i$, as desired.
%The edge set $\mcE$ of the graph $G = (\mcX, \mcE)$ is comprised of
%of those edges $(i,j)$ for which $s_{ij} > 0$ (note that $s_{ij} = s_{ji})$
%and since $P$ is irreducible it follows that $G$ is connected.

We may thus turn to the 
%more interesting 
case when not all the $P(u)$ are the same.

Suppose first that $|\mcX| = 2$, and write $\mcX = \{1,2\}$. 
Then, for any $u \in \mcU$, $P(u)$ is irreducible and
reversible iff we have both $p_{12}(u) > 0$ and $p_{21}(u) > 0$
(the corresponding stationary distribution is
$\left[ \begin{array}{c} \frac{p_{21}(u)}{p_{12}(u) + p_{21}(u)} 
\frac{p_{12}(u)}{p_{12}(u) + p_{21}(u)} \end{array} \right]$).
It can be checked that any collection $(P(u): u \in \mcU)$ where each
$P(u)$ is irreducible and reversible, together with a reward function
$r: \mcX \times \mcU \to \mbbR$, defines 
%a reversible Markov decision problem
an RMDP
(because $P(\mu)$ will then be irreducible and reversible for each
$\mu \in \mcM$). 
We can then define 
$P^\parrzero := \left[ \begin{array}{cc} 0 & 1 \\1 & 0 \end{array} \right]$, 
with $\rho(1,u) = p_{12}(u)$, and
$\rho(2,u) = p_{21}(u)$,
thus establishing the main claim of the theorem in
this case.
Note that the graph associated to this %reversible Markov decision problem 
RMDP
is biconnected.
%We can then take $s_{12} = s_{21} = 1$
%(so that $s_1 = s_2 = 1$ and $S = 2$), $\rho(1,u) = p_{12}(u)$, and
%$\rho(2,u) = p_{21}(u)$, thus establishing the main claim of the theorem in
%this case, even though the graph associated to this reversible Markov decision
%problem is not $2$-connected.

%In general, whatever the cardinality of $\mcX$, 

Having dealt with the case $|\mcX| = 2$, we may henceforth assume that
$|\mcX| \ge 3$.
%as needed.
%We will now make use for the first time of the assumption that $G$ is $2$-connected.
Fix $i \in \mcX$. By the assumption 
%of $2$-connectedness of $G$, 
that $G$ is biconnected 
there must
exist distinct $j, k \in \mcX$ such that $(i,j) \in \mcE$ and $(i,k) \in \mcE$.
This means that for all $u \in \mcU$ we have $p_{ij}(u) > 0$ and $p_{ik}(u) > 0$.
We claim that the ratio $\frac{p_{ij}(u)}{p_{ik}(u)}$ does not depend on $u$. 
To see this, let $u, v \in \mcU$ be distinct and pick some $\tilu \in \mcU$ 
(which could be either $u$ or $v$ if desired)
and consider the two stationary deterministic Markov strategies 
$\barmu^\parra$ and $\barmu^\parrb$ given as in \eqref{eq:twostrats}.
Write $p^\parra_{lm}$ for $p_{lm}(\barmu^\parra)$ for $l,m \in \mcX$,
%and $a \in \{0,1\}$,
and 
%write 
$\pi^\parra$ for the stationary distribution of $P(\barmu^\parra)$;
similarly for $\barmu^\parrb$.
%for $a \in \{0,1\}$.
By the assumption that $G$ is biconnected, there is a path in $G$ from $j$ to $k$ that 
does not touch $i$, i.e. one can find a sequence $(l_0, l_1, \ldots, l_R)$ of vertices of $G$,
where $R \ge 1$, with $l_0 = j$, $l_R = k$, $l_r \neq i$ for $0 \le r \le R$, and such that 
$(l_r, l_{r+1}) \in \mcE$ for $0 \le r \le R-1$. %Then, for each $a \in \{0,1\}$, 
Then
we have the 
equations
\begin{equation}        \label{eq:balanceoutside}
\pi^\parra_{l_r} p^\parra_{l_r l_{r+1}} = \pi^\parra_{l_{r+1}} p^\parra_{l_{r+1} l_r}
\mbox{ and }
\pi^\parrb_{l_r} p^\parrb_{l_r l_{r+1}} = \pi^\parrb_{l_{r+1}} p^\parrb_{l_{r+1} l_r},
\end{equation}
for all $0 \le r \le R-1$ (these follow from the reversibility of $P(\barmu^\parra)$
and $P(\barmu^\parrb)$ respectively). 
Since for all $l, m \in \mcX$ with $l \neq i$ and $m \neq i$ we have
$p^\parra_{lm} = p^\parrb_{lm} = p_{lm}(\tilu)$, we can conclude from 
the equations in \eqref{eq:balanceoutside} that 
\begin{equation}        \label{eq:sameratio}
\frac{\pi^\parra_j}{\pi^\parra_k} = \frac{\pi^\parrb_j}{\pi^\parrb_k}.
\end{equation}
But the reversibility of $P(\barmu^\parra)$ and $P(\barmu^\parrb)$ also gives us the equations
\begin{eqnarray*}
\pi^\parra_i p_{ij}(u) &=& \pi^\parra_j p_{ji}(\tilu),\\
\pi^\parra_i p_{ik}(u) &=& \pi^\parra_k p_{ki}(\tilu),\\
\pi^\parrb_i p_{ij}(v) &=& \pi^\parrb_j p_{ji}(\tilu),\\
\pi^\parrb_i p_{ik}(v) &=& \pi^\parrb_k p_{ki}(\tilu).
\end{eqnarray*}
Dividing the first of these by the second (on each side) and the third of these by the fourth
and comparing the resulting equations, using \eqref{eq:sameratio} we conclude that 
$\frac{p_{ij}(u)}{p_{ik}(u)}$ equals $\frac{p_{ij}(v)}{p_{ik}(v)}$. Since $u, v \in \mcU$, $u \neq v$, 
were arbitrarily chosen, we conclude that $\frac{p_{ij}(u)}{p_{ik}(u)}$ does not depend on $u$, as claimed.

Now, for each $i \in \mcX$, pick 
an arbitrary 
$u_i \in \mcU$ 
%such that $p_{ii}(u_i) \le p_{ii}(u)$ for all $u \in \mcU$.
(all of these could be the same action, if one wishes).
Having made such a choice, define $\barmu \in \bar{\mcM}$ by $\barmu(i) = u_i$ for all $i \in \mcX$.
Since $P(\barmu)$ is irreducible and reversible, we have the equations
$\pi_i(\barmu) p_{ij}(u_i) = \pi_j(\barmu) p_{ji}(u_j)$ for all distinct $i,j \in \mcX$, where
$(\pi_i(\barmu): i \in \mcX)$ denotes the stationary distribution of $P(\barmu)$ as usual. 
For $i \neq j$, 
define $p_{ij}^\parrzero := \frac{p_{ij}(u_i)}{1 - p_{ii}(u_i)}$,
%define $s_{ij} := \pi_i(\barmu) p_{ij}(u_i)$,
and let $p_{ii}^\parrzero = 0$ for all $i \in \mcX$. 
%and let $s_{ii} = 0$ for all $i \in \mcX$. 
Note that $p_{ij}^\parrzero > 0$
%Note that $s_{ij} > 0$ 
iff $(i,j) \in \mcE$, where
$G = (\mcX, \mcE)$ is the 
canonical graph
%simple connected graph 
%without self-loops 
%that has already been 
associated to this
%reversible Markov decision problem.
RMDP.
The resulting matrix $P^\parrzero$ based
on $\mcX$ is an irreducible 
TPM
%transition probability matrix 
with zero diagonal entries,
and it is reversible because its stationary distribution is $(K \pi_i(\barmu) (1 - p_{ii}(u_i)): i \in \mcX)$, where $K$ is the proportionality constant.
%Further, we have $s_{ij} = s_{ji}$, and
%$s_i := \sum_{j \neq i} s_{ij} = \pi_i(\barmu) (1 - p_{ii}(u_i))$ for all $i \in \mcX$.
We can now set $\rho(i,u) = 1 - p_{ii}(u)$ for all $(i,u) \in \mcX \times \mcU$. 
Indeed, 
%since we have already proved that for each $i$ and each $j,k$ distinct from $i$ such that 
%$p_{ij}(u) > 0$ and $p_{ik}(u) > 0$ for some (and hence all) $u \in \mcU$, 
we have already proved that 
the $(p_{ij}(u): j \neq i)$ for $u \in \mcU$
are proportional (for fixed $i \in \mcX$),
and so
we will have
$p_{ij}(u) = p_{ij}(u_i) \frac{1 - p_{ii}(u)}{1 - p_{ij}(u_i)}$ for all $(i,j) \in \mcE$,
%namely $\frac{p_{ij}(u)}{1 - p_{ii}(u)}$ will equal $\frac{s_{ij}}{s_i}$ for all $(i,j) \in \mcE$ and
%all $u \in \mcU$ and so, with $\rho(i,u)$ as defined, 
%and so we will have 
which gives
$\rho(i,u) p_{ij}^\parrzero = p_{ij}(u)$
%$\rho(i,u) \frac{s_{ij}}{s_i} = p_{ij}(u)$
for all $u \in \mcU$ and all distinct $i,j \in \mcX$. Note that we have $\rho(i,u) \in (0,1]$ for all
$(i,u)$, as required.

This concludes the proof of the theorem.
\epf

\section{Proof of Lemma \ref{lem:treecase}}        \label{app:treecase}

For all $\mu \in \mcM$ we have 
$p_{ij}(\mu) > 0$ iff $(i,j) \in \mcE$,
and so $P(\mu)$ is an irreducible 
TPM
%transition probability matrix 
on $\mcX$.
For $i \in \mcX$ and $k \neq i$ define
%$p_{\vec{ki}}(\mu)$ 
$p_{k \to i}(\mu)$ 
to be $p_{kj}(\mu)$,
where $j \in \mcX$ is defined as
the vertex adjacent to $k$ on the unique
path from $k$ to $i$ in the tree.
It can be checked that the stationary distribution of $P(\mu)$ is proportional to
$(\prod_{k \neq i} p_{k \to i}(\mu): i \in \mcX)$ and so $P(\mu)$ is reversible. This concludes the proof.
\epf

\section{Proof of Theorem \ref{thm:maintheorem}}        \label{app:maintheorem}

If there is only one block then we are in biconnected case covered in Theorem \ref{thm:biconnected}, where we have already proved that the structure of the RMDP must be consistent with the type described in Example \ref{ex:genexample}. We may therefore assume that there are at least two blocks, and so every block has at least one articulation point. For each block $B \in \mcB$ an argument similar to that in Theorem \ref{thm:biconnected} shows that for each $i \in B$ the $(p_{ij}(u): j \in B)$ as $u$ ranges over $\mcU$ are all proportional. 
We can therefore find a 
TPM
%transition probability matrix 
$P^\parrzero(B) = \left[ \begin{array}{c} p_{ij}^\parrzero(B) \end{array} \right]$ on $B$, with zero diagonal entries, such that 
$p_{ij}(u) = \left( \sum_{k \in B} p_{ik}(u) \right) p_{ij}^\parrzero$ for all distinct $i,j \in B$. 
Since $p_{ij}^\parrzero(B) > 0$ iff
$(i,j)$ is an edge in $B$ (viewed as a subgraph), and since $B$ is connected, we see that $P^\parrzero(B)$ is irreducible.
Define $\rho(i,u)$ to be
$\sum_{j \in B} p_{ij}(u)$ for $i \in \ringB$
(if any) and, for each articulation point $a \in B$, define $\rho(a,u)$ to be 
$\sum_{j \in \mcX, j \neq a} p_{aj}(u)$
(this quantity does not depend on which $B$ containing $a$ is being considered), and 
define $\nu_a(u,B)$ to be 
$\frac{\sum_{j \in B} p_{ij}(u)}{\rho(a,u)}$.
Note that the $\nu_a(u,B)$ are strictly positive and $\sum_{B \ni a} \nu_a(u,B) = 1$,
as required. Also note that 
$\rho(i,u) \in (0,1]$ for all $(i,u) \in \mcX \times \mcU$. 

It remains to show that each $P^\parrzero(B)$ is reversible.
Pick any $u \in \mcU$. 
Let $(\pi_i(u): i \in \mcX)$
denote the stationary distribution of $P(u)$.
Fix $B \in \mcB$.
By the reversibility of $P(u)$
we have $\pi_i(u) p_{ij}(u) = \pi_j(u) p_{ji}(u)$ for all $i,j \in B$. 
It follows that 
\[
\left( \pi_i(u) \sum_{k \in B} p_{ik}(u) \right) p_{ij}^\parrzero = 
\left( \pi_j(u) \sum_{k \in B} p_{jk}(u) \right) p_{ji}^\parrzero 
\]
This means that if
$(\psi_i(B): i \in B)$ 
denotes the stationary distribution of 
$P^\parrzero(B)$ then it is proportional
to $(\pi_i(u) \sum_{k \in B} p_{ik}(u): i \in B)$ and thus that 
$\psi_i(B) p_{ij}^\parrzero(B) = 
\psi_j (B) p_{ji}^\parrzero(B)$ for all
$i, j \in B$, which establishes that
$P^\parrzero(B)$ is reversible. This concludes the proof.
\epf

\section{Proof of Theorem \ref{thm:policyiterationbic}}        \label{app:policyiterationbic}

Let $(\pi_i^\parrzero: i \in \mcX)$ denote the stationary distribution of $P^\parrzero$, and 
recall 
%from the proof of Theorem \ref{thm:biconnected} 
that for any $\mu \in \mcM$ the stationary distribution of $P(\mu)$ is given by
$(K(\mu) \frac{\pi_i^\parrzero}{\rho(i,\mu(i))}: i \in \mcX)$, where $K(\mu) := \left( \sum_i \frac{\pi_i^\parrzero}{\rho(i,\mu(i))} \right)^{-1}$ is the normalizing constant.
Since
\[
\beta(\mu) = \sum_i r(i, \mu(i)) \pi_i(\mu) = K(\mu) \sum_i r(i, \mu(i)) \frac{\pi_i^\parrzero}{\rho(i,\mu(i))},
\]
we get
\[
\sum_i \frac{r(i, \mu(i)) - \beta(\mu)}{\rho(i,\mu(i))}
\pi_i^\parrzero= 0.
\]
Thus $\beta(\mu)$ can be characterized as
\[
\beta(\mu) = \sup \{ \beta \in \mbbR: 
\sum_i \frac{r(i, \mu(i)) - \beta}{\rho(i,\mu(i))} \pi_i^\parrzero \ge 0\}.
\]

If we can find $i \in \mcX$ for which 
equation \eqref{eq:badbic} holds, then pick $u \in \mcU$ achieving the argmax on the RHS of equation \eqref{eq:badbic} and let $\barmu^{(k+1)}(i) = u$ and $\barmu^{(k+1)}(j) = \barmu^\parrk(j)$ for all $j \neq i$, as in the simplified policy iteration algorithm.
We then have
\begin{eqnarray*}
\sum_j \frac{r(j, \barmu^{(k+1)}(j)) - \beta(\barmu^\parrk)}{\rho(j,\barmu^{(k+1)}(j))} \pi_j^\parrzero
&=&
\sum_{j \neq i} \frac{r(j, \barmu^{(k)}(j)) - \beta(\barmu^\parrk)}{\rho(j,\barmu^{(k)}(j))} \pi_j^\parrzero
+ 
\frac{r(i, u) - \beta(\barmu^\parrk)}{\rho(i,u)} \pi_i^\parrzero\\
&>&
\sum_{j} \frac{r(j, \barmu^{(k)}(j)) - \beta(\barmu^\parrk)}{\rho(j,\barmu^{(k)}(j))} \pi_j^\parrzero
= 0.
\end{eqnarray*}
It follows that 
$\beta(\barmu^{(k+1)}) > 
\beta(\barmu^\parrk)$, which concludes the proof.

%Since the summation on the right hand side of equation \eqref{eq:DPeqn} can be restricted to $j$ distinct from $i$, substituting for $p_{ij}(v)$ with 
%$\rho(i,v) p_{ij}^\parrzero$ this equation reads
%\[
%0 = \max_v \left(
%\frac{r(i,v) - \beta(\mu)}{\rho(i,v)}
%+ \sum_{j \neq i} p_{ij}^\parrzero(h_j(\mu) - h_i(\mu))
%\right),
%\]
%which must hold for all $i \in \mcX$ for any optimal $\mu \in \mcM$. But then this means that we must have
%\[
%0 = \max_v 
%\frac{r(i,v) - \beta(\mu)}{\rho(i,v)},
%\]
%for all $i \in \mcX$ for any optimal $\mu \in \mcM$.

%Let $(\pi_i^\parrzero: i \in \mcX)$ denote the stationary distribution of $P^\parrzero$, and 
%recall 
%from the proof of Theorem \ref{thm:biconnected} 
%that for any $\mu \in \mcM$ the stationary distribution of $P(\mu)$ is given by
%$(K(\mu) \frac{\pi_i^\parrzero}{\rho(i,\mu(i))}: i \in \mcX)$, where $K(\mu) := \left( \sum_i \frac{\pi_i^\parrzero}{\rho(i,\mu(i))} \right)^{-1}$ is the normalizing constant.
%We then have
%\begin{eqnarray*}
%\beta(\barmu^{(k+1)}) &=& \sum_j r(j, \barmu^{(k+1)(j)}) \pi_j(\barmu^{(k+1)})\\
%&=& \sum_j r(j, \barmu^{(k+1)(j)}) K(\barmu^{(k+1)}) \frac{\pi_j^\parrzero}{\rho(j,\barmu^{(k+1)}(j))}
%\end{eqnarray*}

\section{Proof of Theorem \ref{thm:policyiterationgen}}        \label{app:policyiterationgen}

Define $p_{ij}^\parrzero(\barmu|_\mcA)$ to be
$\nu_i(\barmu(i), B) p_{ij}^\parrzero(B)$ for $i, j \in B$, $i \neq j$, 
for each $B \in \mcB$. Here we recall that we defined
$\nu_i(u,B) = 1$ for all $i \in \ringB$ and $u \in \mcU$, and so we realize that 
 $p_{ij}^\parrzero(\barmu|_\mcA)$ depends only on the restriction of
 $\barmu$ to the articulation nodes, which is indicated
by the notation $\barmu|_\mcA$. It is straighforward to check that
the $p_{ij}^\parrzero(\barmu|_\mcA)$ define a 
TPM
%transition probability matrix 
on $\mcX$. Let $(\pi_i^\parrzero(\barmu|_\mcA): i \in \mcX)$ denote the stationary distribution associated to this TPM.
It is straightforward to check that the stationary distribution of 
$P(\barmu)$ is proportional to 
$(\frac{\pi_i^\parrzero(\barmu|_\mcA)}{\rho(i,\barmu(i))}:  i \in \mcX)$.
%Further, because $k \in \ringB$ for some $B \in \mcB$ (so 
%$\nu_k(u,B) = 1$ for all $u \in \mcU$, we can check that the 
%stationary distribution of 
%$P(\barmu^\prime)$ is proportional to 
%$(\frac{\pi_i^\parrzero(\barmu|_\mcA)}{\rho(i,\barmu^\prime(i))}:  i \in \mcX)$.
Further, we can check that for all $\bareta \in \bar{\mcM}$ such that
$\bareta|_\mcA = \barmu|_\mcA$
 the stationary distribution of 
$P(\barmu)$ is proportional to 
$(\frac{\pi_i^\parrzero(\barmu|_\mcA)}{\rho(i,\bareta(i))}:  i \in \mcX)$.

From this, as in the proof of Theorem \ref{thm:policyiterationbic}, 
we can check that for all  $\bareta \in \bar{\mcM}$ such that
$\bareta|_\mcA = \barmu|_\mcA$ we have the characterization
\[
\beta(\bareta) = \sup \{ \beta \in \mbbR: 
\sum_i \frac{r(i, \bareta(i)) - \beta}{\rho(i,\bareta(i))} 
\pi_i^\parrzero(\barmu|_\mcA) \ge 0\}.
\]
The rest of the proof then follows as in the proof of Theorem 
\ref{thm:policyiterationbic}, allowing us to conclude the desired strict inequality.
\epf

\iffalse
\section{Rough Notes}

{\em Note 1}: The multivariate Gaussian with covariance matrix
$\pi(i) Z(i,j)$ where $Z$ is the fundamental matrix
ought to show up from the Gaussian free field in the construction of Janson (i.e. requiring only that sums around loops are zero, which leaves the mean value free)
by conditioning on $\sum_i \pi(i) \eta(i) = 0$. Check on this.\\

{\em Note 2}: Need to figure out how to write the matrix
$\pi_{\mu}(i) Z_{\mu}(i,j)$ in terms of the matrix
$\pi(i) Z(i,j)$, where $\mu: \mcX \to \mcU$ is a deterministic control strategy.
\fi

\iffalse
$\pi^\parrj$

$\pi^\chevzero$

$\pi^\curlthree$

$\pi^\sqrrnine$

$\pi^\parreight$
\fi


\begin{thebibliography}{99}

\bibitem{AldousFill}
David Aldous and James Allen Fill.
{\em Reversible Markov Chains and Random Walks on Graphs}.
Unfinished monograph (2002). Recompiled version, 2014.
https://www.stat.berkeley.edu/users/aldous/RWG/book.pdf

\bibitem{ABFGM93}
Aristotle Arapostathis, Vivek S. Borkar, Emmanuel Fern\'{a}ndez-Gaucherand, Mrinal K. Ghosh, and Steven I. Marcus.
``Discrete-time controlled Markov processes with average cost criterion: A survey."
{\em SIAM Journal on Control and Optimization}, Vol. 31, No. 2, pp. 282 -344, 1993.

\bibitem{CP2013}
Randy
Cogill and Cheng
Peng.
``Reversible Markov decision processes with an average-reward criterion."
{\em SIAM Journal on Control and Optimization}, vol. 51, No. 1, 
pp. 402 -418, 2013.

\bibitem{Diestel}
Reinhard Diestel.
{\em Graph Theory}.
Fifth edition, Springer, 2017.

\bibitem{LeJanSF}
Yves Le Jan.
{\em Markov Path, Loops and Fields}.
\'{E}cole d'\'{e}t\'{e} de Probabilit\'{e}s
de Saint-Flour, XXXVIII, 2008. Springer, 2011.

\bibitem{LPbook}
Russell Lyons and Yuval Peres.
{\em Probability on Trees and Networks}.
Cambridge University Press, 2016.

\bibitem{Janson}
Svante Janson.
{\em Gaussian Hilbert Spaces}.
Cambridge Tracts in Mathematics, Vol. 129. Cambridge University Press,
1997.

\bibitem{SznitmanETH}
Alain-Sol
Sznitman.
{\em Topics in Occupation Times and Gaussian Free Fields}.
Notes of the course ``Special topics in probability"
at ETH Zurich, Spring 2011. 


\end{thebibliography}
\end{document}